\documentclass[11pt,a4paper]{article}
\usepackage[utf8]{inputenc}
\usepackage{amsmath}
\usepackage{amsfonts}  
\usepackage{amsmath}
\usepackage{amssymb}
\usepackage{amsthm}
\usepackage{mathrsfs}
\usepackage{enumerate}
\usepackage{mathabx}
\usepackage{hyperref}
\usepackage[scr=rsfso]{mathalfa}
\usepackage{fancyhdr}
\usepackage{etoolbox}

\usepackage[a4paper, right=3.2cm, left=3.2cm, top=4.0cm, bottom=3.9cm]{geometry}

\usepackage{titlesec}
\titleformat{\section}{\large\bfseries}{\thesection}{1em}{}
\numberwithin{equation}{section}

\newtheoremstyle{noenddot}
{\topsep}
{\topsep}
{\itshape}
{-4pt}
{\bfseries}
{}
{ } 
{\thmname{#1}\thmnumber{ #2}\thmnote{ \normalfont(#3)}}

\theoremstyle{noenddot}
\newtheorem{theorem}{}[section]

\newcommand{\thm}{\textnormal{\textbf{Theorem.~}}}
\newcommand{\proposition}{\textnormal{\textbf{Proposition.~}}}

\newcommand{\lemma}{\textnormal{\textbf{Lemma.~}}}
\newcommand{\definition}{\textnormal{\textbf{Definition.~}}}

\newcommand{\1}{1\!\!1}

\newcommand{\overbar}[1]{\mkern 3mu\overline{\mkern-3mu#1\mkern-0.5mu}}

\title{\vspace{-5ex}\large\textbf{STRICT IRREDUCIBILITY OF MARKOV CHAINS AND ERGODICITY OF SKEW PRODUCTS}}
\author{\normalsize\textsc{Pablo Lummerzheim, Felix Pogorzelski and Elias Zimmermann}}
\date{\vspace{-5ex}}

\begin{document}
	
	\maketitle
	
	\begin{abstract}
		\textsc{Abstract}. We consider a family of measure preserving transformations, which act on a common probability space and are chosen at random by a stationary ergodic Markov chain. This setting defines an instance of a random dynamical system, which may be described in terms of a step skew product. In many contexts it is desirable to know whether ergodicity of the family implies ergodicity of the skew product. Introducing the notion of strict irreducibility for Markov kernels we shall characterize the class of Markov chains for which the aforementioned implication holds true. We thereby extend a sufficient condition of Bufetov for the case of finite state Markov chains to general state spaces and show that it is in fact also necessary. As an application we obtain an explicit description of the limit in ergodic theorems for a suitable class of random transformations.
	\end{abstract} 
	
	\section{Introduction}
	
	Let $(E,\mathscr{E})$ and $(X,\mathscr{A})$ be measurable spaces. Fixing a probability space $(\Omega,\mathscr{C},\nu)$ we shall consider a (measurable) family of transformations $(T_{y})_{y \in E}$ on $X$, which are chosen randomly by a stationary and ergodic stochastic process $(\xi_{n})_{n=0}^{\infty}$ consisting of random variables $\xi_{n}\colon \Omega \to E$. This setting describes an instance of a random dynamical system. Systems of the above type play an important role in various areas of ergodic theory and related fields, see e.\@ g. \cite{Arn98}, \cite{Fur02} or \cite{Kif86}. A natural object of study in this context are random ergodic theorems, i.\@ e. statements about almost sure (or other types of) convergence of the averages 
			\[A_{n}f := \frac{1}{n}\sum_{i=0}^{n-1}f \circ T_{\xi}^{i}\]
	of a function $f \in L^{1}(\mu)$ along the random compositions $\smash{T_{\xi}^{n}} := \smash{T_{\xi_{n-1}} \circ \ldots \circ T_{\xi_{0}}}$ for $n \to \infty$. The study of random ergodic theorems goes back to Pitt and Ulam-von Neuman, see \cite{Pit42} and \cite{UvN45}, and has been an important research topic ever since. A theorem establishing mean and pointwise convergence of random averages for a large class of random dynamical systems is due to Kin, cf. \cite{Kin72}. However, there is no explicit description of the limit function. In what follows we will assume that the family of transformations $(T_{y})_{y \in E}$ is measure preserving and ergodic with respect to a common probability measure $\mu$ on $X$. This setting arises for instance naturally in the study of measure preserving group actions, cf. \cite{BK12} and \cite{Fur02}. Under these assumptions one would expect the limit of the above averages to be the integral $\int f~\!d\mu$. The classical random ergodic theorem of this type for i.i.d. processes is due to Kakutani and Ryll-Nardzewski, cf. \cite{Kak51} and \cite{Ryll55}. However, beyond the i.i.d. setting the identification of the limit in the random ergodic theorem is an open problem in many situations.
	
	In what follows we shall suppose that the random variables $\xi_{n}$ are given in their canonical form, i.\@ e. $(\Omega,\mathscr{C})$ is the product space $(E^{\mathbb{N}_{0}},\mathscr{E}^{\otimes\mathbb{N}_{0}})$ and $\xi_{n}$ is the projection to the $n$-th coordinate. In this setting stationarity (ergodicity) of the process $(\xi_{n})_{n=0}^{\infty}$ amounts to the assumption that the measure $\nu$ is invariant (ergodic) with respect to the left shift $S$ on $\Omega$. An important tool for the study of random ergodic theorems is the investigation of the step skew product $T$ on $\Omega \times X$ arising from $S$ and $(T_{y})_{y \in E}$. In the above setting it is easy to see that $T$ preserves the product measure $\nu \otimes \mu$. Observing that the above averages are just the ordinary ergodic averages of the function $\1 \otimes f$ with respect to the skew product $T$, we obtain that $A_{n}f$ converges $\nu \otimes \mu$-almost surely to some limit function $\smash{\overbar{f}} \in L^{1}(\nu \otimes \mu)$ by Birkhoff's ergodic theorem. If the skew product $T$ is moreover ergodic, it is not difficult to see that $\smash{\overbar{f}}$ is constant a.\@ s. and equals the integral $\int f~\!d\mu$. Thus we have reduced the above problem to the question whether every skew extension $T$ of $\xi$ with an ergodic family of transformations $(T_{y})_{y \in E}$ is ergodic. 
	
	Kakutani showed in \cite{Kak51} that if  $\xi$ is an i.i.d. process and the transformations are invertible then the aforementioned implication holds true. The assumption of invertibility was later removed by Ryll-Nardzewski, see \cite{Ryll55}. For Markov chains with finite state space Bufetov introduced in \cite{Buf00a}, \cite{Buf00b} and \cite{Buf01} the concept of strict irreducibility and showed that it provides a sufficient criterion for the above property to hold. 
	
	The main goal of this paper is to generalize Bufetov's result to Markov chains with arbitrary state space. To this end we introduce the notion of strict irreducibility for general Markov kernels, cf. Definition \ref{DefStrictIrreducibility}. A key role in the definition is played by the concept of a deterministic set, which is introduced in the same definition. In our main result, Theorem \ref{main}, we shall show that for a strictly irreducible Markov chain every skew extension with an ergodic family of measure preserving transformations is ergodic. This contains the results of Kakutani and of Bufetov as special cases. To this end we use a result of Kowalski, cf. \cite{Kow15}, characterizing eigenfunctions of the Perron-Frobenius operator associated to a skew product of the above type. Moreover we will show that if a Markov chain satisfies the above property, then it is necessarily strictly irreducible. As an application we obtain a substantial generalization of Kakutani's classical random ergodic theorem for i.i.d. processes to the setting of Markov chains. 
	
	\medskip
	
	The results of this work also appeared in the diploma thesis of one of the authors (P.\@ L.), supervised by the other two authors (F.\@ P.) and (E.\@ Z.). This diploma project was pivotal to extend the setting under consideration from Markov shifts over countable state spaces to Markov shifts over arbitrary measurable state spaces.
	
	\medskip
	
	The paper is organized as follows: Section \ref{SectionPrelim} recalls some preliminaries on Markov operators and dynamical systems. In Section \ref{SectionStrictIrreducibility} we fix the setting and give a definition of strict irreducibility for Markov kernels, cf. Definition \ref{DefStrictIrreducibility}. Section \ref{SectionSkewProduct} is devoted to the proof of our main result, Theorem \ref{main}. Finally, in Section \ref{SectionErgodicTheorem}, we derive two versions of the random ergodic theorem for transformations, which are selected by Markov chains, see Theorems \ref{random} and \ref{random2}.
	
	\section{Preliminaries} \label{SectionPrelim}
	
	In this section we collect some basic facts about Markov operators and dynamical systems. Let $(Y,\mathscr{B},\eta)$ be a probability space. By a \textit{measure preserving transformation} we mean a measurable map $T\colon Y \to Y$ such that $\eta(T^{-1}D) = \eta(D)$ for all sets $D \in \mathscr{B}$. The triple $(Y,\eta,T)$ is called a \textit{measure preserving dynamical system (MDS)}. A set $D \in \mathscr{B}$ is \textit{invariant} under $T$, or \textit{$T$-invariant} for short, if it satisfies $\eta(D \Delta T^{-1}D) = 0$. The collection of $T$-invariant sets forms a $\sigma$-algebra, which shall be denoted by $\Sigma_{T}$. The transformation $T$ is said to be \textit{ergodic} if $\Sigma_{T}$ is trivial, i.\@ e. it contains only sets of trivial measure.
	
	Given a bounded linear operator $B\colon L^{p}(\eta) \to L^{p}(\eta)$, where $p \in [1,\infty)$, we shall denote by $B'$ the \textit{dual operator}, i.\@ e. the unique bounded linear operator $B'\colon L^{q}(\eta) \to L^{q}(\eta)$ such that 
	\[\int_{Y}Bg \cdot h~d\eta = \int_{Y}g \cdot B'h~d\eta\]
	for all $g \in L^{p}(\eta)$ and $h \in L^{q}(\eta)$, where $q$ denotes the dual exponent of $p$. For ${p = 2}$ we shall furthermore denote by $B^{*}$ the \textit{adjoint operator}, i.\@ e.\@ the unique bounded linear operator $B^{*}\colon L^{2}(\eta) \to L^{2}(\eta)$ satisfying $\langle Bg,h\rangle =  \langle g,B^{*}h\rangle$ for all $g \in L^{2}(\eta)$ and $h \in L^{2}(\eta)$, where $\langle \cdot,\cdot\rangle$ denotes the standard inner product of $L^{2}(\eta)$. Let $T$ be a measure preserving transformation on $(Y,\eta)$. For any $p \in [1,\infty]$ the \textit{Koopman operator} $\smash{\widehat{T}}\colon L^{p}(\eta) \to L^{p}(\eta)$ corresponding to $T$ is defined by $\smash{\widehat{T}} g := g \circ T$ for $g \in L^{p}(\eta)$. It is easy to see that $\smash{\widehat{T}}$ is linear and isometric. By a \textit{$T$-invariant} function we mean a fixed function of $\smash{\widehat{T}}$. It is not difficult to verify that a function $f \in L^{p}(\eta)$ is $T$-invariant if and only if it is $\Sigma_{T}$-measurable. Considering the Koopman operator $\smash{\widehat{T}}$ on $L^{\infty}(\eta)$ there is a unique bounded linear operator $\mathcal{L}_{T}\colon L^{1}(\eta) \to L^{1}(\eta)$ satisfying $\mathcal{L}_{T}' = \smash{\widehat{T}}$, see \cite[§1.3]{Aar97}. $\mathcal{L}_{T}$ is called the \textit{Perron-Frobenius operator} corresponding to $T$. 
	
	Given a measurable function $f$ we shall write $f \vargeq 0$ to indicate that $\eta$-almost everywhere $f$ is real-valued and non-negative. A linear operator $M \colon L^{2}(\eta) \to L^{2}(\eta)$ is called \textit{positive} if $f \vargeq 0$ implies $Mf \vargeq 0$ for all $f \in L^{2}(\eta)$. Positive operators are always bounded, cf. \cite[Lemma 7.5]{EFHN15}. By a \textit{Markov operator} we mean a positive operator $M$ satisfying $M\1 = \1$ and
	\[\int_{Y} Mf~\!d\eta = \int_{Y} f~\!d\eta\]
	for all $f \in L^{2}(\eta)$. Markov operators are contractions with operator norm equal to $1$ and the class of Markov operators is closed under composition and taking adjoints, cf. \cite[Section 13.1]{EFHN15}. We call a Markov operator $M$ \textit{irreducible} if the only $M$-invariant functions are the constants. This notion can be generalized to positive operators on arbitrary Banach lattices, cf. \cite[Chapter 7]{EFHN15}. On the level of sets it means that $M\1_{D} = \1_{D}$ implies $\eta(D) \in \{0,1\}$ for all sets $D \in \mathscr{B}$, where one may equivalently replace the condition $M\1_{D} = \1_{D}$ by the weaker condition $M\1_{D} \varleq \1_{D}$. 
	
	\section{Strict irreducibility of Markov kernels} \label{SectionStrictIrreducibility}
	
	In this section we develop the relevant background on Markov chains. Moreover, we introduce the notion of strict irreducibility for Markov kernels, which will play the key role throughout what follows. To motivate our definition, let us start by recalling the notion of strict irreducibility for Markov chains with finite state space as introduced by Bufetov in \cite{Buf00a}, \cite{Buf00b} and \cite{Buf01}. 
	
	Let $\Pi$ be the transition matrix and $m$ be the initial probability vector corresponding to a Markov chain $\xi$ with finite state space $E$. The Markov chain $\xi$ is stationary if and only if $\pi$ is a fixed vector of $\Pi^{T}$. In what follows we will assume that $\pi$ has only positive entries. Since in the stationary setting states with probability zero are never reached this is no restriction of generality. Recall that a row-stochastic matrix $\Pi$ is called \textit{irreducible} if there is some $n \in \mathbb{N}$ such that $\Pi + \Pi^{2} + \ldots + \Pi^{n}$ has only positive entries. A stationary Markov chain $\xi$ is ergodic if and only if $\Pi$ is irreducible. In this case the Markov chain itself is also called \textit{irreducible}. As we shall see below it is not difficult to show that for an irreducible matrix $\Pi$ the matrix $\Pi^{T}\Pi$ is irreducible if and only if the matrix $\Pi~\!\Pi^{T}$ is irreducible. Following Bufetov we call an irreducible matrix $\Pi$ satisfying the latter condition \textit{strictly irreducible}. Accordingly a stationary irreducible Markov chain $\xi$ is called \textit{strictly irreducible} if the transition matrix $\Pi$ is strictly irreducible. 
	
	It is well known that irreducibility of a row stochastic matrix $\Pi$ can be characterized in terms of the connectedness of the associated transition graph. A characterization of the same type can also be obtained for strict irreducibility of $\Pi$. Denoting by $\pi(i,j)$ the $(i,j)$-entry of $\Pi$ consider the relation $\sim$ on $E$, which takes place between two states $i$ and $j$ if there is some state $k \in E$ such that $\pi(k,i) > 0$ and $\pi(k,j) > 0$. Obviously, this defines a symmetric relation. We may turn $\sim$ into an equivalence relation by considering its transitive closure $\simeq$, i.\@ e. we have $i \simeq j$ if and only if $i = j$ or $i \sim k_{1}, k_{1} \sim k_{2},\ldots,k_{r} \sim j$ for some sequence of states $k_{1},\ldots,k_{r} \in E$. It is not difficult to see that the graph $(E,\sim)$ is connected if and only if the equivalence relation $\simeq$ has only one equivalence class. Moreover, a straightforward computation shows that both conditions are equivalent to the irreducibility of the matrix $\Pi^{T}\Pi$. An analogous characterization of the irreducibility of the matrix $\Pi~\!\Pi^{T}$ can be obtained by considering the dual relation $\backsim$, which takes place between two states $i,j \in E$ if there is a state $k \in E$ such that $\pi(i,k) > 0$ and $\pi(j,k) > 0$. However, if $\Pi$ is irreducible, it is not difficult to see that both characterizations describe equivalent conditions, which shows that in this case $\Pi^{T}\Pi$ is irreducible if and only if $\Pi~\!\Pi^{T}$ is irreducible.
	
	In \cite{Buf00b} and \cite{Buf01} strict irreducibility is defined using the latter condition. In \cite{Buf00a} the terminology of strongly and $*$-strongly connected Markov measures is used to refer to conditions, which are equivalent to the two combinatorial characterizations of strict irreducibility given before. By the foregoing discussion, all of these notions describe equivalent concepts.
	
	Strictly irreducible finite state Markov chains have proved to be a powerful tool for the study of actions of certain word hyperbolic groups. The approach is to extend classical random walk techniques as used e.\@ g.\@ in \cite{Ose65},\cite{BN15}, \cite{AL05} and \cite{Kai03} to suitable non-backtracking random walks in order to prove ergodic theorems or to study mixing properties of boundary actions. For instance, based on the concept of strict irreducibility, Bufetov developed an operator theoretic approach to ergodic theorems for group actions. The latter was successfully applied for improving and generalizing the well-known spherical ergodic theorems for free groups due to Nevo-Stein and Grigorchuk, cf. \cite{NS94}, \cite{Gri86} and \cite{Gri99}. It also allowed for an extension of these results to a much larger class of hyperbolic groups, including certain right angled Artin groups and Fuchsian groups, see \cite{BK12}, \cite{BK12a}, \cite{BKK12}, \cite{BS11} and \cite{BKS23}. Recently the concept of strict irreducibility has also been linked to mixing properties of boundary actions. The authors have been informed that in the upcoming preprint \cite{TZ} of Tserunyan and Zomback a characterization of weak mixing for the boundary action of finitely generated free groups and semigroups with respect to Markov measures was obtained, showing in particular that it is equivalent to strict irreducibility of the transition matrix.
	
	Using the combinatorial characterization from above the notion of strict irreducibility can canonically be extended to Markov chains with countably infinite state space. Indeed, Grigorenko introduced an equivalent condition for Markov chains with (possibly infinite) discrete state space and showed that such Markov chains satisfy a Hewitt-Savage type theorem in the sense that their symmetric $\sigma$-algebra is trivial, see \cite{Gri86}. On the other hand, it is a priori not clear how to generalize the above concepts to Markov chains with arbitrary measurable state space. Obviously the above conditions cannot directly be translated into this context. However, we can obtain a further characterization of strict irreducibility, which allows for such a generalization.
	
	To this end let us call a set $B \subseteq E$ \textit{deterministic} if $\pi(i,B) \in \{0,1\}$ for any $i \in E$. It is not hard to see that the deterministic sets are precisely the unions of connected components of the graph $(E,\sim)$. Thus we obtain that $\Pi$ is strictly irreducible if and only if $\emptyset$ and $E$ are the only deterministic sets. The notion of a deterministic set can also be defined for general Markov kernels as we shall see below, cf. Definition \ref{DefStrictIrreducibility}, and this opens a way to extend the notion of strict irreducibility to Markov chains with general state space. In Proposition \ref{StrictIrreducibility} below we will see that the same notion can also be expressed in terms of the Markov operator associated to the Markov chain, which provides us with a condition analogous to the definition of strict irreducibility in terms of transition matrices. Such a condition has been studied by Bezhaeva and Oseledets in \cite{BO96} in order to describe the symmetric $\sigma$-algebra of Harris-Markov chains. 
	
	\medskip
	
	Before we go into the details of the approach just described we have to develop some background on Markov kernels and Markov operators. To this end let $(E,\mathscr{E})$ be a measurable space. A map $\pi\colon E \times \mathscr{E} \to [0,1]$ is called a \textit{Markov kernel} if for all $y \in E$ the map $\pi(y,\cdot)$ is a probability measure, which will occasionally be denoted by $\pi_{y}$ in the following, and for all $B \in \mathscr{E}$ the map $\pi(\cdot,B)$ is measurable. The product of two Markov kernels $\pi$ and $\kappa$ is given by
	\[\pi\kappa(y,B) := \int_{E}\kappa(z,B)~d\pi_{y}(z)\]
	for $y \in E$ and $B \in \mathscr{E}$ and defines again a Markov kernel. We shall call a probability measure $m$ on $\mathscr{E}$ \textit{invariant} under a Markov kernel $\pi$, or \textit{$\pi$-invariant} for short, if it satisfies 
	\[m(B) = \int_{E}\pi(y,B)~dm(y)\]
	for all $B \in \mathscr{E}$. It is not difficult to see that a probability measure $m$, which is invariant under two Markov kernels $\pi$ and $\kappa$, is also invariant under the product kernel $\pi\kappa$. 
	
	Given a Markov kernel $\pi$ together with a $\pi$-invariant probability measure $m$ we may define a bounded linear operator $P\colon L^{2}(m) \to L^{2}(m)$, which is a well known example of a Markov operator, by setting
	\[P f(y) := \int_{E}f(z)~d\pi_{y}(z)\]
	for $f \in L^{2}(m)$ and $y \in E$, see \cite[Section 1.6]{DMPS18}. The product of Markov kernels is compatible with the composition of the respective Markov operators in the following sense: If $P$ and $Q$ are the Markov operators corresponding to Markov kernels $\pi$ and $\kappa$ with joint invariant probability measure $m$, then the Markov operator corresponding to the product kernel $\pi\kappa$ is given by the composition $PQ$ of $P$ and $Q$. 
	
	Let us introduce the following notation: Given a set $B \in \mathscr{E}$ together with an exponent $n \in \mathbb{N}$ we shall denote by $\smash{U_{B}^{~\!\!n}}$ the set of all $y \in E$ satisfying $\pi^{n}(y,B) > 0$. Noting that $\smash{U_{B}^{~\!\!n}}$ is just the preimage of the set $(0,1]$ under the map $\pi^{n}(\cdot,B)$ we may conclude that $\smash{U_{B}^{~\!\!n}}$ is measurable for all $n \in \mathbb{N}$ and so is the set
	\[U_{B} := \bigcup_{n=1}^{\infty}U_{B}^{~\!\! n}.\]
	We will derive some properties of $U_{B}$ that will be used repeatedly in what follows. First of all, we may note that for all $y \in U_{B}^{~\!\!c}$ and all $n \in \mathbb{N}$ we have
	\begin{align*} 
		\int_{E}\pi^{n}(z,B)~d\pi_{y}(z) = \pi^{n+1}(y,B) = 0.
	\end{align*}
	Since the functions $\pi^{n}(\cdot,B)$ are non-negative, this implies that $\pi(y,U_{B}^{~\!\!n}) = 0$ for all $n \in  \mathbb{N}$ and thus $\pi(y,U_{B}) = 0$ for every $y \in U_{B}^{~\!\!c}$. Combining this fact with the $\pi$-invariance of $m$ we obtain furthermore
	\begin{align*} \label{eq7}
		\int_{U_{B}}\pi(y,U_{B})~dm(y) = \int_{E}\pi(y,U_{B})~dm(y) = m(U_{B})
	\end{align*}
	and thus $\pi(y,U_{B}) = 1$ for $m$-almost all $y \in U_{B}$. Finally, using again the $\pi$-invariance of $m$, we get
	\begin{align*} 
		m(B) = \int_{E} \pi(y,B)~dm(y) = \int_{U^{~\!\!1}_{B}}\pi(y,B)~dm(y).
	\end{align*} 
	Therefore, if $B$ has positive measure, then $U_{B}^{1}$ and thus $U_{B}$ have positive measure too. 
	
    \begin{theorem} \definition \textup{We shall call a Markov kernel $\pi$ \textit{irreducible} with respect to a $\pi$-invariant probability measure $m$ if for every set $B \in \mathscr{E}$ with $m(B) > 0$ $m$-almost all states $y \in E$ admit some $n \in \mathbb{N}$ (which may depend on $y$) such that $\pi^{n}(y,B) > 0$.}
    \end{theorem}
    In terms of the notation introduced above the latter just means that $m(U_{B}) = 1$. The following proposition connects this notion of irreducibility to the notion of irreducibility for Markov operators. 
    
	\begin{theorem} \label{prop1} \proposition
		A Markov kernel $\pi$ is irreducible with respect to a $\pi$-invariant proba\-bility measure $m$ if and only if the Markov operator $P$ corresponding to $\pi$ and $m$ is irreducible.
	\end{theorem}
	
	\textit{Proof:} Assume that $P$ is not irreducible. Then there is some set $B \in \mathscr{E}$ with $m(B) \in (0,1)$ such that $P\1_{B} = \1_{B}$. This implies that for $m$-almost all $y \in B^{c}$ and every $n \in \mathbb{N}$ we have  
	\[\pi^{n}(y,B) = P^{n}\1_{B}(y) = \1_{B}(y) = 0.\]
	Since $B$ and $B^{c}$ have both positive measure, the latter contradicts the irreducibility of $\pi$ with respect to $m$. 
	
	To show the converse direction assume that $P$ is irreducible and fix a set $B \in \mathscr{E}$ with $m(B) > 0$. We have to show that $m(U_{B}) = 1$. Since $m$ is $\pi$-invariant, we know that $m(U_{B}) > 0$. Furthermore we have $P\1_{U_{B}}(y) = \pi(y,U_{B}) = 0$ for $m$-almost all $y \in U_{B}^{~\!\!c}$ and $P\1_{U_{B}}(y) = \pi(y,U_{B}) = 1$ for $m$-almost all $y \in U_{B}$, which together gives $P\1_{U_{B}} = \1_{U_{B}}$. By the irreducibility of $P$ this implies that $U_{B}$ has trivial measure. Since we already know that $m(U_{B}) > 0$, we obtain $m(U_{B}) = 1$. 
	\hfill $\Box$ \\
	
	The above proposition shows that the notion of irreducibility of Markov kernels introduced above is tightly related to the notion of $m$-irreducibility in the theory of Markov chains, cf. \cite[Section 2.2]{Num84} and \cite[Section 9.2]{DMPS18}. However, the latter is slightly stronger in the sense that it requires for any $B \in \mathscr{E}$ with $m(B) > 0$ that \textit{all} $y \in E$ admit some $n \in \mathbb{N}$ (which may depend on $y$) such that $\pi^{n}(y,B) > 0$. It will turn out that from a dynamical point of view the weaker notion of irreducibility is more convenient. It is for instance equivalent to the ergodicity of the Markov shift corresponding to $\pi$ and $m$ as we shall see in Proposition \ref{ergodicity} below.  
	
	\medskip
	
	As indicated above we shall denote by $(\Omega,\mathscr{C})$ the product space $(E^{\mathbb{N}_{0}},\mathscr{E}^{\otimes\mathbb{N}_{0}})$. Consider a Markov kernel $\pi$ and a probability measure $m$ as before. Then, as a consequence of Carathéodory's extension theorem, there is a unique probability measure $\nu$ on $\mathscr{C}$ satisfying
	\[\nu\big(B_{0} \times \cdots \times B_{k-1} \times \Omega\big) = \int_{B_{0}}\ldots\int_{B_{k-1}}~d\pi_{y_{k-2}}(y_{k-1})\ldots d\pi_{y_{0}}(y_{1})dm(y_{0})\]
	for every $k \in \mathbb{N}$ and all sets $B_{0},\ldots,B_{k-1} \in \mathscr{E}$, cf. \cite[Theorem 3.1.2]{DMPS18}. The measure $\nu$ is called the \textit{Markov measure} corresponding to $\pi$ and $m$. The \textit{shift map} $S\colon \Omega \to \Omega$ given by
	\[S(\omega_{0}\omega_{1}\ldots) := \omega_{1}\omega_{2}\ldots\]
	for $\omega \in \Omega$ defines a $\mathscr{C}$-measurable map. We shall call the Markov measure $\nu$ \textit{stationary} if it is preserved under $S$. It is easy to check that $\nu$ is stationary if and only if the measure $m$ is $\pi$-invariant. In this case we obtain an MDS $(\Omega,\nu,S)$, which we shall call a \textit{Markov shift}. 
	
	For any probability measure $m$ we may consider the trivial kernel $\tau$, which is given by $\tau(y,B) := m(B)$ for $y \in E$ and $B \in \mathscr{E}$. Then $m$ is $\tau$-invariant. In this case the corresponding Markov measure coincides with the product measure $m^{\otimes\mathbb{N}_{0}}$ and the corresponding Markov chain defines an i.i.d. process. We will refer to the respective Markov shift as a \textit{Bernoulli shift}. The next proposition connects the ergodicity of the Markov shift $S$ to the irreducibility of the kernel $\pi$ with respect to the measure $m$. Noting that trivial kernels are irreducible this contains the well known fact that Bernoulli shifts are ergodic as a special case.
	
	\begin{theorem} \proposition \label{ergodicity}
		Let $\pi$ be a Markov kernel with invariant probability measure $m$. Then the respective Markov shift $S$ is ergodic if and only if $\pi$ is irreducible with respect to $m$.
	\end{theorem}
	
	\textit{Proof:} Recall that a set $B \in \mathscr{E}$ is called absorbing if we have $\pi(y,B) = 1$ for all $y \in B$. It is well known that the Markov shift $S$ corresponding to $\pi$ and $m$ is ergodic if and only if $m(B) \in \{0,1\}$ for every absorbing set $B \in \mathscr{E}$, see \cite[Theorem 5.2.11]{DMPS18}. 
	
	Assume that $S$ is not ergodic. Then there is an absorbing set $B \in \mathscr{E}$ with ${m(B) \in (0,1)}$. We claim that $\pi^{n}(y,B^{c}) = 0$ for all $y \in B$ and $n \in \mathbb{N}$. This can be shown inductively as follows: For $n = 1$ the claim follows from the definition of an absorbing set. Now for $n > 1$ we have 
	\begin{align*}
		\pi^{n}(y,B^{c}) = \int_{B}\pi(z,B^{c})~d\pi_{y}^{n-1}(z) + \int_{B^{c}}\pi(z,B^{c})~d\pi_{y}^{n-1}(z)
	\end{align*}
	for all $y \in B$. The function integrated in the first term vanishes on $B$ by the fact that $B$ is absorbing. Furthermore, by the induction hypothesis, we have $\pi^{n-1}(y,B^{c}) = 0$ for all $y \in B$, so we integrate over a null set in the second term. Therefore both terms are zero, which gives $\pi^{n}(y,B^{c}) = 0$ for all $y \in B$. Since $B$ and $B^{c}$ have positive measure by assumption, we obtain that $\pi$ is not irreducible with respect to $m$. 
	
	Conversely assume that $\pi$ is not irreducible with respect to $m$. Then there is a set $B \in \mathscr{E}$ with $m(B) > 0$ such that $U_{B}^{~\!\!c}$ has positive measure. Since $m$ is $\pi$-invariant, $U_{B}$ has also positive measure, so we get $m(U_{B}^{~\!\!c}) \in (0,1)$. Recall that we have $\pi(y,U_{B}) = 0$ and thus $\pi(y,U^{~\!\!c}_{B}) = 1$ for all $y \in U_{B}^{~\!\!c}$, so $U_{B}^{~\!\!c}$ is a non trivial absorbing set. This shows that $S$ is not ergodic. \hfill $\Box$ \\  
	
	We shall now introduce the central notion of this paper, namely the concept of strict irreducibility for general Markov kernels, in terms of deterministic sets.
    
    \begin{theorem} \definition \label{DefStrictIrreducibility}
    \textup{Given a Markov kernel $\pi$ admitting a $\pi$-invariant probability measure $m$ we will call a set $B \in \mathscr{E}$ \textit{deterministic} if for $m$-almost all $y \in E$ we have either $\pi(y,B) = 0$ or $\pi(y,B) = 1$. Based on this we shall say that $\pi$ is \textit{strictly irreducible} with respect to $m$ if for all deterministic sets $B \in \mathscr{E}$ we have $m(B) \in \{0,1\}$.}
    \end{theorem}
    
    Note that in the case of a discrete state space the particular choice of $m$ does not matter as long as $m$ gives positive mass to all points (which we may always assume in the stationary setting), so in this case the definition coincides with the definition of a deterministic set given in the introduction of this section and we obtain the notion of strict irreducibility for discrete state Markov chains as a special case.
	
	As we shall see below, the condition given above can also be expressed in terms of the Markov operator $P$ corresponding to $\pi$ and $m$. To this end we will rely on the following general property of irreducible Markov operators.
	
	\begin{theorem} \label{lemma} \lemma
		Let $(Y,\eta)$ be a probability space and $M\colon L^{2}(\eta) \to L^{2}(\eta)$ be an irreducible Markov operator. Consider functions $g,h \in L^{2}(\eta)$ with $g,h \vargeq 0$ such that $g + h = \1$ and $\langle Mg,h\rangle = 0$. Then either $g = 0$ or $h = 0$. 
	\end{theorem}
	
	\textit{Proof:} Let $D$ denote the preimage of the interval $(0,\infty)$ under $h$. Then $D$ is obviously measurable. Furthermore, since $Mg \vargeq 0$, $\langle Mg,h\rangle = 0$ implies that $\1_{D}Mg = 0$. On the other hand we have 
	\[Mg + Mh = M(g+h) = M\1 = \1,\]
	so we obtain 
	\[\1_{D}Mh = \1_{D}Mg + \1_{D}Mh = \1_{D}(Mg + Mh) = \1_{D}.\]
	Both together yields
	\begin{align*}
		\int_{D}h~\!d\eta = \int_{Y}h~\!d\eta = \int_{Y}Mh~\!d\eta \vargeq \int_{D}Mh~\!d\eta = \int_{D} \1~\!d\eta = \eta(D).
	\end{align*}
	The assumptions on $g$ and $h$ imply that $h \varleq \1$, so we may conclude that $h = \1_{D}$. This gives $g = \1_{D^{c}}$ and thus 
	\[\1_{D}M\1_{D^{c}} = \1_{D}Mg = 0.\] 
	Using that $Mh \vargeq 0$ we obtain
	\[M\1_{D^{c}} = Mg = \1 - Mh \varleq \1.\]
	Both together implies 
	\[M\1_{D^{c}} = \1_{D}M\1_{D^{c}} + \1_{D^{c}}M\1_{D^{c}} \varleq \1_{D^{c}}.\]
	By the irreducibility of $M$ it follows that $D^{c}$ and thus $D$ has trivial measure. Accordingly we obtain that either $g=0$ or $h=0$. \hfill $\Box$ \\
	
	The above lemma allows us to characterize the property of strict irreducibility in terms of the irreducibility of Markov operators as follows. 
	
	\begin{theorem} \label{StrictIrreducibility} \proposition
		Let $\pi$ be a Markov kernel and $m$ be a $\pi$-invariant probability measure. Let $P$ denote the Markov operator corresponding to $\pi$ and $m$. Then $P^{*}P$ is irreducible if and only if $PP^{*}$ is irreducible. Both is the case if and only if $\pi$ is strictly irreducible with respect to $m$.
	\end{theorem}
	
	\textit{Proof:} To show the first equivalence assume that $P^{*}P$ is irreducible and fix a set $B \in \mathscr{E}$ with $PP^{*}\1_{B} = \1_{B}$. Then we have 
	\[\langle P^{*}\1_{B},P^{*}\1_{B^{c}}\rangle = \langle PP^{*}\1_{B},\1_{B^{c}}\rangle = \langle \1_{B},\1_{B^{c}}\rangle = 0.\]
	Since $P^{*}\1_{B} \vargeq 0$ and $P^{*}\1_{B^{c}} \vargeq 0$ this implies $P^{*}\1_{B} P^{*}\1_{B^{c}} = 0$. Accordingly, using that $P^{*}\1_{B} + P^{*}\1_{B^{c}} = P^{*}\1 = \1$, we obtain $P^{*}\1_{B} = \1_{C}$ and $P^{*}\1_{B^{c}} = \1_{C^{c}}$ for some set $C \in \mathscr{E}$. It is easy to see that we also have $PP^{*}\1_{B^{c}} = \1_{B^{c}}$ by the fact that $PP^{*}$ is a Markov operator. Thus we get
	\begin{align*}
		\langle P^{*}P\1_{C},\1_{C^{c}}\rangle &= \langle P^{*}PP^{*}\1_{B},P^{*}\1_{B^{c}}\rangle = \langle PP^{*}\1_{B},PP^{*}\1_{B^{c}}\rangle = \langle \1_{B},\1_{B^{c}} \rangle = 0.
	\end{align*}
	Since $P^{*}P$ is irreducible by assumption, we may apply Lemma \ref{lemma} to obtain that ${m(C) \in \{0,1\}}$. If $m(C) = 1$ this gives $P^{*}\1_{B} = \1$. However, since $P^{*}$ is bounded with $\|P^{*}\| = \|P\| = 1$, this is only possible if $m(B) = 1$. In the case $m(C) = 0$ we may use the same argument to obtain that $m(B^{c}) = 1$. In sum we see that $B$ has trivial measure. This shows that $PP^{*}$ is irreducible. A symmetric argument gives the reverse implication.
	
	To show the second equivalence assume again that $P^{*}P$ is irreducible. Let $B \in \mathscr{E}$ be a deterministic set. Then we have $\pi(y,B)\pi(y,B^{c}) = 0$ for $m$-almost all $y \in E$ and thus
	\[\langle P^{*}P\1_{B},\1_{B^{c}}\rangle = \langle P\1_{B},P\1_{B^{c}}\rangle = \int_{E}\pi(y,B)\pi(y,B^{c})~dm(y) = 0,\]
	so by Lemma \ref{lemma} we obtain $m(B) \in \{0,1\}$. This shows that $\pi$ is strictly irreducible with respect to $m$. To show the reverse direction assume that $\pi$ is strictly irreducible with respect to $m$ and fix a set $B \in \mathscr{E}$ with $P^{*}P\1_{B} = \1_{B}$. Then we have
	\begin{align*}
		\smash{\int_{E}}\pi(y,B)\pi(y,B^{c})~dm(y) &=  \langle P\1_{B},P\1_{B^{c}}\rangle = \langle P^{*}P\1_{B},\1_{B^{c}}\rangle = \langle \1_{B},\1_{B^{c}}\rangle = 0.
	\end{align*}
	This implies
	\[\pi(y,B)\pi(y,B^{c}) = 0\]
	and thus $\pi(y,B) \in \{0,1\}$ for $m$-almost all $y \in E$. Accordingly $B$ is deterministic, so we get $m(B) \in \{0,1\}$. This proves that $P^{*}P$ is irreducible. \hfill $\Box$ \\
	
	Note that in the proof of the first statement of Proposition \ref{StrictIrreducibility} we did not use any specific properties of the operator $P$. In fact the equivalence claimed in this statement holds for arbitrary Markov operators on $L^{2}$-spaces.
	
	\medskip
	
	Let us remark that in the case that $(E,\mathscr{E})$ is a standard Borel space the adjoint operator $P^{*}$ of $P$ admits a more concrete description, which gives rise to a further characterization of strict irreducibility. This is due to the fact that in this setting there exists a Markov kernel $\pi^{*}$ such that
	\begin{align*}
		\int_{E}\int_{E}f(y,z)~d\pi_{y}(z)dm(y) = \int_{E}\int_{E}f(y,z)~d\pi^{*}_{z}(y)dm(z)
	\end{align*}
	for all $f \in L^{2}(m)$. Such a kernel is called a \textit{reverse kernel} of $\pi$. Indeed, the kernels $\pi$ and $\pi^{*}$ coincide with the regular conditional distributions of the unique measure $\eta$ on $\mathscr{E} \otimes \mathscr{E}$ satisfying 
	\[\eta(B_{0} \times B_{1}) = \int_{B_{0}}\pi(y,B_{1})~dm(y)\]
	for all $B_{0},B_{1} \in \mathscr{E}$, which exist by the disintegration theorem for measures whenever $(E,\mathscr{E})$ is a standard Borel space, cf. \cite[Chapter 17]{Kec95}. The disintegration theorem guarantees moreover that the reverse kernel $\pi^{*}$ is unique up to $m$-null sets, i.\@ e. for any other reverse kernel $\kappa$ one has $\kappa_{y} = \pi^{*}_{y}$ for $m$-almost all $y \in E$. 
	
	It is readily seen that $m$ is also $\pi^{*}$-invariant. Denote by $P^{*}$ the Markov operator corresponding to $\pi^{*}$ and $m$. By the essential uniqueness of the reverse kernel $P^{*}$ is unique as a Markov operator on $L^{2}(m)$.  As an immediate consequence of the defining equation of the reverse kernel one obtains then 
	\[\langle Pg,h\rangle = \langle g,P^{*}h\rangle\]
	for all $g,h \in L^{2}(m)$. This implies that $P^{*}$ coincides indeed with the adjoint operator of $P$. By Proposition \ref{StrictIrreducibility} and the foregoing results this shows that $\pi^{*}\pi$ is irreducible with respect to $m$ if and only if $\pi\pi^{*}$ is irreducible with respect to $m$ and both is equivalent to the strict irreducibility of $\pi$ with respect to $m$.
	
	\medskip
	
	The next proposition shows that strict irreducibility is indeed a stronger property than irreducibility. In fact, using Lemma \ref{lemma}, it is not difficult to give an abstract argument for this fact on the level of Markov operators. However, the proof of the following proposition, which operates on the level of Markov kernels, might be more instructive in the present setting.
	
	\begin{theorem} \label{StrictlyIrreducibleImpliesIrreducible} \proposition
		A Markov kernel $\pi$ which is strictly irreducible with respect to an invariant probability measure $m$ is also irreducible with respect to $m$.
	\end{theorem}
	
	\textit{Proof:} Assume that $\pi$ is not irreducible with respect to $m$. Then there is some set $B \in \mathscr{E}$ with $m(B) > 0$ such that $m(U_{B}) < 1$. However, as noted several times already, we have $m(U_{B}) > 0$ by the $\pi$-invariance of $m$ as well as $\pi(y,U_{B}) = 0$ for all $y \in U_{B}^{~\!\!c}$ and $\pi(y,U_{B}) = 1$ for $m$-almost all $y \in U_{B}$. This shows that $U_{B}$ is a non-trivial deterministic set. Thus $\pi$ is not strictly irreducible with respect to $m$. \hfill $\Box$

	\section{Ergodicity of step skew products} \label{SectionSkewProduct}
	
	This section is devoted to the proof of our main result, Theorem \ref{main}, which states that a stationary Markov chain is strictly irreducible if and only if every skew extension with an ergodic family of measure preserving transformations is ergodic. This extends Bufetov's result for Markov chains with finite state space to general state spaces and shows that the condition of strict irreducibility is not only sufficient but also necessary for the above property to hold. Recalling that i.i.d. processes are special cases of strictly irreducible Markov chains it provides also a generalization of the results of Kakutani and Ryll-Nardzewski.
	
	Several arguments in Bufetov's original proof rely on the assumption of a discrete state space and it is not obvious how to generalize them to the general setting. We will therefore follow a different proof strategy. In order to study the invariant functions of the step skew product we will make use of the spectral theory of the associated Perron-Frobenius operator. More precisely, we will utilize a result of Kowalski, cf.
	\cite{Kow15}, providing a convenient representation of the eigenfunctions of this operator. Using the characterization of strict irreducibility in terms of Markov operators obtained in the last section this will allow us to obtain a particularly simple description of the invariant functions of a step skew product over a strictly irreducible Markov chain, cf. Theorem \ref{representation2}. This representation will be the central ingredient in the proof of Theorem \ref{main}.
	
	Let $(E,\mathscr{E})$ be a measurable space. We call a family $(T_{y})_{y \in E}$ of transformations $T_{y}$ on a further measurable space $(X,\mathscr{A})$ \textit{measurable} if the map $(y,x) \mapsto T_{y}(x)$ is measurable with respect to $\mathscr{E} \otimes \mathscr{A}$ and $\mathscr{A}$. Considering a Markov shift $(\Omega,\nu,S)$ with state space $E$ such a family gives rise to a measurable skew product $T$ on $\Omega \otimes X$ defined by
	\[T(\omega,x) := (S\omega,T_{\omega_{0}}x)\]
	for $(\omega,x) \in \Omega \times X$. Skew products of the above form are often called \textit{step skew products}. The mapping $(n,\omega,x) \mapsto T_{\omega}^{n}(x)$, where $T_{\omega}^{n}$ denotes the random composition $T_{\omega_{n-1}} \circ \ldots \circ T_{\omega_{0}}$ for $n \in \mathbb{N}_{0}$, defines a \textit{(measurable) random dynamical system} in the sense of \cite{Arn98}. In what follows we will assume that the transformations $(T_{y})_{y \in E}$ preserve a common probability measure $\mu$ on $X$. It is easily verified that in this case $T$ preserves the product measure $\nu \otimes \mu$.
	
	It will turn out that the invariant functions of $T$ are of a particular simple form. This is a consequence of a result of Kowalski concerning the eigenfunctions of the Perron-Frobenius operator associated to step skew product as above, which may be stated as follows. Let $\mathcal{L}_{T}$ denote the Perron-Frobenius operator of $T$ and consider a function $\varphi \in L^{1}(\nu \otimes \mu)$ with $\mathcal{L}_{T}\varphi = \lambda\varphi$ for some $\lambda \in \mathbb{C}$ with $|\lambda| = 1$. Then by \cite[Theorem 3.1]{Kow15} there is a function $\widehat{\varphi} \in L^{1}(m \otimes \mu)$ satisfying 
	\begin{align*}
		\varphi(\omega,x) = \widehat{\varphi}(\omega_{0},x)
	\end{align*}
	for $\nu \otimes \mu$-almost all $(\omega,x) \in \Omega \times X$. In other words, in the $\Omega$-argument the eigenfunctions of $\mathcal{L}_{T}$ depend only on the first symbol. This property transfers to $T$-invariant functions as follows.
	
	\begin{theorem} \proposition \label{representation1}
		Let $(\Omega,\nu,S)$ be a Markov shift with state space $E$ and $(T_{y})_{y \in E}$ be a measurable family of measure preserving transformations on a probability space $(X,\mu)$. Let $T$ be the arising step skew product and $\varphi \in L^{1}(\nu \otimes \mu)$ be a $T$-invariant function. Then there is a function $\widehat{\varphi} \in L^{1}(m \otimes \mu)$ such that for $\nu \otimes \mu$-almost all $(\omega,x) \in \Omega \times X$ we have
		\[\varphi(\omega,x) = \widehat{\varphi}(\omega_{0},x).\]
	\end{theorem}
	
	\textit{Proof:} Consider the Koopman operator $\widehat{T}$ on $L^{\infty}(\nu \otimes \mu)$ corresponding to $T$ and let $\varphi \in L^{1}(\nu \otimes \mu)$ be a $T$-invariant function. Then by the $T$-invariance of the measure $\nu \otimes \mu$ we obtain 
	\begin{align*}
		\int_{\Omega \times X} {\cal L}_{T}\varphi\cdot\psi~d\nu \otimes \mu &= 	\int_{\Omega \times X} \varphi\cdot\widehat{T}\psi~d\nu \otimes \mu  \\
		&= \int_{\Omega \times X} \varphi\circ T\cdot \psi\circ T~d\nu\otimes\mu 
		\\ &= \int_{\Omega \times X} \varphi\cdot \psi~d\nu\otimes\mu 
	\end{align*}
	for all $\psi \in L^{\infty}(\nu \otimes \mu)$. This shows that ${\cal L}_{T}\varphi = \varphi$, so the claimed representation of $\varphi$ follows from Kowalski's theorem. \hfill $\Box$ \\
	
	Let $\Sigma$ denote the $\sigma$-algebra of sets $A \in \mathscr{A}$ satisfying $\mu(A \Delta T_{y}^{-1}A) = 0$ for $m$-almost all $y \in E$. We shall call the family $(T_{y})_{y \in E}$ \textit{ergodic} if $\Sigma$ is trivial. It is not difficult to see that the ergodicity of the skew product $T$ implies the ergodicity of the family $(T_{y})_{y \in E}$ as well as the ergodicity of the Markov shift $S$. However, the converse does not hold in general. To see this consider a $3$-point set $\overbar{X} := \{1,2,3\}$ equipped with the uniform probability measure $\overbar{\mu}$ giving every point equal weight. Consider furthermore the Markov chain with state space $\{0,1\}$ corresponding to the transition matrix
	\[\Pi = \left(\begin{array}{cc}0 & 1 \\ 1 & 0 \end{array}\right)\]
	and the ($\Pi$-invariant) initial probability vector $m =  \smash{\big(\frac{1}{2},\frac{1}{2}\big)^{T}}$. Then $\Pi$ is irreducible, so the arising Markov shift $(\{0,1\}^{\mathbb{N}_{0}},\overbar{\nu},S)$ is ergodic. Let $T_{0}$ be a fixpoint free permutation of $\overbar{X}$ and set $T_{1} := T_{0}^{-1}$. Then $T_{0}$ and $T_{1}$ are measure preserving and ergodic with respect to $\overbar{\mu}$, in particular the family $\{T_{0},T_{1}\}$ is ergodic. To see that the arising step skew product $T$ is not ergodic observe that the measure $\overbar{\nu}$ is concentrated on the sequences $\omega^{0} := 010101...$ and $\omega^{1} := 101010...$ with equal weight $\frac{1}{2}$. Moreover, it is not difficult to check that the set $\big\{(\omega^{0},1),(\omega^{1},T_{0}1)\big\}$ is $T$-invariant. However, this set has measure $\frac{1}{3}$. So $T$ admits a non-trivial invariant set and can thus not be ergodic.
	
	In Theorem \ref{main} below we shall see that the converse implication holds true if the Markov shift arises from a strictly irreducible Markov kernel (in fact, it is easy to see that the matrix $\Pi$ in the example above is not strictly irreducible). In order to obtain such a result we have to improve Proposition \ref{representation1}. To this end let us introduce the following notation. Fixing $y \in E$ let $\delta_{y}$ be the Dirac measure concentrated in $y$ and let $\pi$ be any Markov kernel. We denote by $\nu_{y}$ the Markov measure arising from $\delta_{y}$ and $\pi$, which by definition takes the form 
	\begin{align} \label{eq3}
		\nu_{y}(B_{0} \times \dots \times B_{k-1} \times \Omega) = \1_{B_{0}}(y)\int_{B_{1}}\ldots\int_{B_{k-1}} d\pi_{z_{k-2}}(z_{k-1})\ldots d\pi_{y}(z_{1})
	\end{align}
	for all sets $B_{0},\ldots,B_{k-1} \in \mathscr{E}$ and every $k \in \mathbb{N}$. Based on the above identity a standard approximation argument shows that for any $k \in \mathbb{N}$ and any bounded measurable function $h\colon E^{k} \to \mathbb{R}$ we have
	\begin{align} \label{eq2}
		\int_{\Omega}h(\omega_{0},\ldots,\omega_{k-1})~\!d\nu_{y}(\omega) =\int_{E}\ldots\int_{E}h(y,z_{1},\ldots,z_{k-1})~\!d\pi_{z_{k-2}}(z_{k-1})\ldots d\pi_{y}(z_{1}).
	\end{align}
	Furthermore, one can show that the family of measures $(\nu_{y})_{y \in E}$ provides a disintegration of the measure $\nu$ with respect to $m$, cf. \cite[Proposition 3.1.3]{DMPS18}. This implies that for all sets $C \in \mathscr{C}$ we have 
	\begin{align} \label{eq0}
		\nu(C) = \int_{E} \nu_{y}(C)~\!dm(y).
	\end{align}
	
	We are now able to improve the representation of $T$-invariant functions given in Proposition \ref{representation1} for step skew products over Markov shifts in case of a strictly irreducible Markov kernel. It will turn out that in this situation the $T$-invariant functions do not depend on the $\Omega$-argument at all.
	
	\begin{theorem} \label{representation2} \thm
		Let $(\Omega,\nu,S)$ be a Markov shift with state space $E$ and $(T_{y})_{y \in E}$ be a measurable family of measure preserving transformations on a probability space $(X,\mu)$. Let $T$ denote the respective step skew product. Assume that the Markov measure $\nu$ arises from a Markov kernel $\pi$ and a $\pi$-invariant probability measure $m$ such that $\pi$ is strictly irreducible with respect to $m$.  Then for every $T$-invariant function $\varphi \in L^{1}(\nu \otimes \mu)$ there is a function $f \in L^{1}(\mu)$ such that for $\nu \otimes \mu$-almost all $(\omega,x)$ we have 
			\[\varphi(\omega,x) = f(x).\]
	\end{theorem}
	
	\textit{Proof:} It will suffice to prove the above statement for all characteristic functions of $T$-invariant sets. This is due to the fact that a function is $T$-invariant if and only if it is $\Sigma_{T}$-measurable, so every $T$-invariant function is the pointwise limit of linear combinations of such characteristic functions.
	
	Let $\varphi := \1_{D}$ be the characteristic function of a $T$-invariant set $D$ and set $\psi := \1_{D^{c}}$. Then Proposition \ref{representation1} provides us with functions $\smash{\widehat{\varphi}}, \smash{\widehat{\psi}} \in L^{1}(m \otimes \mu)$ such that for $\nu \otimes \mu$-almost all $(\omega,x) \in \Omega \times X$ we have 
	\[\varphi(\omega,x) = \widehat{\varphi}(\omega_{0},x)\]
	and 
	\[\psi(\omega,x) = \widehat{\psi}(\omega_{0},x).\]
	This implies that there is a set $B \in \mathscr{E} \otimes \mathscr{A}$ such that $\widehat{\varphi} = \1_{B}$ and $\smash{\widehat{\psi}} = \1_{B^{c}}$. It remains to show that up to measure zero we have
	$B = E \times A$
	for some set $A \in \mathscr{A}$, which gives the claimed representation of $\varphi$ with $f := \1_{A}$. 
	
	To this end observe that by (\ref{eq0}) a property holds $\nu$-almost surely if and only if it holds $\nu_{y}$-almost surely for $m$-almost all $y \in E$. Therefore the $T$-invariance of $\varphi$ implies that for $m \otimes \mu$-almost all $(y,x) \in E \times X$ we have  
	\[\varphi(\omega,x) = \varphi(S\omega,T_{\omega_{0}}x)\]
	for $\nu_{y}$-almost every $\omega \in \Omega$. Setting $B_{x} := \{y \in E\colon (y,x) \in B\}$ for $x \in X$ we obtain thus
	\begin{align*}
		\widehat{\varphi}(y,x) &= \1_{B_{x}}(y) = \nu_{y}(B_{x} \times \Omega) = \int_{\Omega}\1_{B_{x}}(\omega_{0})~d\nu_{y}(\omega) = \int_{\Omega}\widehat{\varphi}(\omega_{0},x)~d\nu_{y}(\omega) \\ &= \int_{\Omega}\varphi(\omega,x)~d\nu_{y}(\omega) = \int_{\Omega} \varphi(S\omega,T_{\omega_{0}}x)~d\nu_{y}(\omega) \\ &= \int_{\Omega}\widehat{\varphi}(\omega_{1},T_{\omega_{0}}x)~d\nu_{y}(\omega) = \int_{E}\widehat{\varphi}(z,T_{y}x)~d\pi_{y}(z)
	\end{align*}
	for $m \otimes \mu$-almost all $(y,x) \in E \times X$ by (\ref{eq3}) and (\ref{eq2}). A similar argument shows that 
	\[\widehat{\psi}(y,x) = \int_{E}\widehat{\psi}(z,T_{y}x)~d\pi_{y}(z)\]
	for $m \otimes \mu$-almost all $(y,x) \in E \times X$. Setting $\smash{\widehat{\varphi}}^{~\!x} := \smash{\widehat{\varphi}}(\cdot,x)$ and $\smash{\widehat{\psi}}^{~\!x} := \smash{\widehat{\psi}}(\cdot,x)$ for $x \in X$ and noting that $\smash{\widehat{\varphi}}, \smash{\widehat{\psi}} \in L^{2}(m \otimes \mu)$ as well as ${\smash{\widehat{\varphi}}\cdot\smash{\widehat{\psi}} = \1_{B}\cdot\1_{B^{c}} = 0}$ we obtain therefore
	\begin{align*}
		0 &= \langle~ \widehat{\varphi},\widehat{\psi}~\!\rangle = \int_{E}\int_{X}\int_{E}\widehat{\varphi}(z,T_{y}x)~d\pi_{y}(z)\int_{E}\widehat{\psi}(z',T_{y}x)~d\pi_{y}(z')~d\mu(x)dm(y)  \\
		&= \int_{E}\int_{X}\int_{E}\widehat{\varphi}(z,x)~d\pi_{y}(z)\int_{E}\widehat{\psi}(z',x)~d\pi_{y}(z')~d\mu(x)dm(y) \\
		&= \int_{X}\int_{E}P\widehat{\varphi}^{~\!x}(y)P\widehat{\psi}^{~\!x}(y)~dm(y)d\mu(x) = \int_{X}\langle P\widehat{\varphi}^{~\!x},P\widehat{\psi}^{~\!x} \rangle~d\mu(x)
	\end{align*}
	by Fubini's theorem, where we have used that the transformation $T_{y}$ is measure preserving for all $y \in Y$. Since $P$ is a Markov operator and $\smash{\widehat{\varphi}^{~\!x}}$, $\smash{\widehat{\psi}^{~\!x}}$ are non-negative functions, we have $\smash{P\widehat{\varphi}^{~\!x}} \vargeq 0$ and $\smash{P\widehat{\psi}^{~\!x}} \vargeq 0$, so the integrated function in the last integral is non-negative. This implies that for $\mu$-almost all $x \in X$ we get 
	\[\langle P^{*}P\widehat{\varphi}^{~\!x},\widehat{\psi}^{~\!x}\rangle = \langle P\widehat{\varphi}^{~\!x},P\widehat{\psi}^{~\!x}\rangle = 0.\]
	Observing that $\smash{\widehat{\varphi}}^{~\!x} + \smash{\widehat{\psi}}^{~\!x} = \1_{B}(\cdot,x) + \1_{B^{c}}(\cdot,x) = \1$ for $\mu$-almost all $x \in X$ and using that $P^{*}P$ is irreducible by assumption we may apply Lemma \ref{lemma} to obtain that for $\mu$-almost all $x \in X$ we have either $\smash{\widehat{\varphi}}^{~\!x} = 0$ or $\smash{\widehat{\psi}}^{~\!x} = 0$. Thus we may conclude that for every such $x \in X$ we have either 
	\[\1_{B}(y,x) = \widehat{\varphi}(y,x) = 0\]
	for $m$-almost all $y \in E$ or 
	\[\1_{B}(y,x) = 1 - \1_{B^{c}}(y,x) = 1 - \widehat{\psi}(y,x) = 1\]
	for $m$-almost all $y \in E$. This shows that up to a $\nu \otimes \mu$-null set we have $B = E \times A$ for some set $A \in \mathscr{A}$, which finishes the proof.  \hfill $\Box$ \\
	
	Using the above representation of invariant functions for step skew products in case of strictly irreducible Markov kernels we are now able to prove the main theorem of the paper.
	
	\begin{theorem} \label{main} \thm
		Consider a Markov kernel $\pi$ admitting a $\pi$-invariant probability measure $m$ on $E$, let $\nu$ be the respective Markov measure and $S$ denote the corresponding Markov shift. Then the following statements are equivalent:  
		\begin{enumerate}[~~~~~~(i)]
			\item $\pi$ is strictly irreducible with respect to $m$.
			\item For every step skew product $T$ over $S$ arising from a measurable family $(T_{y})_{y \in E}$ of measure preserving transformations on some probability space $(X,\mu)$ we have (up to null sets)
			\[\Sigma_{T} = \Sigma_{S} \otimes \Sigma.\]
			\item $S$ is ergodic and any step skew product $T$ over $S$ arising from an ergodic measurable family $(T_{y})_{y \in E}$ of measure preserving transformations on some probability space $(X,\mu)$ is ergodic. 
		\end{enumerate}
	\end{theorem}
	
	\textit{Proof:}  $(i) \Rightarrow (ii)$: Observing that the measurable rectangles in $\Sigma_{S} \otimes \Sigma$ are $T$-invariant we immediately obtain that $\Sigma_{S} \otimes \Sigma \subseteq \Sigma_{T}$. So it remains to show that up to null sets we have $\Sigma_{T} \subseteq \Sigma_{S} \otimes \Sigma$. To this end fix a set $C \in \Sigma_{T}$ and consider the characteristic function $\varphi := \1_{C}$. Then by Theorem \ref{representation2} we have $\varphi = \1 \otimes f$ for some $f \in L^{1}(\mu)$. Obviously $f = \1_{A}$ for some set $A \in \mathscr{A}$ and thus $C = \Omega \times A$ up to a null set. We have to show that $f$ and thus $A$ is invariant under the family $(T_{y})_{y \in E}$. To this end we compute
	\begin{align*}
		\int_{E}\int_{X}\big| f \circ T_{y}(x) - &f(x)\big|~d\mu(x) dm(y) = \int_{\Omega}\int_{X}\big|f \circ T_{\omega_{0}}(x) - f(x)\big|~d\mu(x)d\nu(\omega) \\
		&= \int_{\Omega}\int_{X} \big|(\1 \otimes f)\circ T(\omega,x) - \1 \otimes f(\omega,x)\big|~d\mu(x)d\nu(\omega) \\
		&= \int_{\Omega \times X}\big|\varphi \circ T - \varphi\big|~d\nu \otimes \mu = 0.
	\end{align*}
	This implies $f \circ T_{y} = f$ for $m$-almost all $y \in E$, which shows the claim.
	
	\medskip
	
	$(ii) \Rightarrow (iii)$ Assume that $S$ is not ergodic. Then by \cite[Theorem 5.2.11]{DMPS18} there is an absorbing set $B \in \mathscr{E}$ with $m(B) \in (0,1)$. Set $C := B \times \Omega$. Then we have $\nu(C) = m(B) \in (0,1)$ and 
	\[\nu(C \cap S^{-1}C^{c}) = \nu(B \times B^{c} \times \Omega) = \int_{B}\pi(y,B^{c})~dm(y) = 0.\]
	Since $S$ is $\nu$-preserving this shows that $C$ is $S$-invariant. Now consider the $2$-point set $\overbar{X} := \{1,2\}$ equipped with the equidistributed probability measure $\overbar{\mu}$ and let $\sigma$ be the non-trivial permutation of $\overbar{X}$. We may define a measurable family $(T_{y})_{y \in E}$ of measure preserving transformations by setting $T_{y} := \text{Id}$ if $y \in B$ and $T_{y} := \sigma$ if $y \in B^{c}$. It is not difficult to see that $\Sigma = \{\emptyset,\overbar{X}\}$, so we get $\Sigma_{S} \otimes \Sigma = \{U \times \overbar{X}\colon U \in \Sigma_{S}\}$. In particular the set $D := C \times \{1\}$ differs from any set in $\Sigma_{S} \otimes \Sigma$ by a set of positive measure. However, using the $S$-invariance of $C$, it is easy to see that $D$ is $T$-invariant and thus lies in $\Sigma_{T}$. This contradicts $(ii)$. Therefore $S$ has to be ergodic. If in addition the family $(T_{y})_{y \in E}$ is ergodic, then $\Sigma_{S}$ and $\Sigma$ are trivial. But by $(ii)$ this implies that $\Sigma_{T}$ is trivial, so $T$ is ergodic. 
	
	\medskip
	
	$(iii) \Rightarrow (i)$ Assume that $\pi$ is not strictly irreducible with respect to $m$. If $\pi$ is not even irreducible with respect to $m$, then $S$ is not ergodic and we are done. So we may assume that $\pi$ is irreducible with respect $m$ and thus that $S$ is ergodic. We shall construct an ergodic measurable family $(T_{y})_{y \in E}$ of measure preserving transformations such that the corresponding skew product $T$ is not ergodic. 
	
	To this end consider again the $2$-point set $\overbar{X} := \{1,2\}$ equipped with the equidistribution $\overbar{\mu}$ and let $\sigma$ denote the non-trivial permutation of $\overbar{X}$. Then $\sigma$ is measure preserving and ergodic with respect to $\overbar{\mu}$. Since $\pi$ is not strictly irreducible with respect to $m$, there is a deterministic set $B \in \mathscr{E}$ with probability $m(B) \in (0,1)$. Consider the sets $B_{1},\ldots,B_{4}$ given by 
	\[B_{1} := \big\{y \in B\colon \pi(y,B) = 1\big\},~B_{2} := \big\{y \in B\colon \pi(y,B) = 0\big\}\]
	and 
	\[B_{3} := \big\{y \in B^{c}\colon \pi(y,B^{c}) = 1\big\},~B_{4} := \big\{y \in B^{c}\colon \pi(y,B^{c}) = 0\big\}.\]
	Since the maps $\pi(\cdot,B)$ and $\pi(\cdot,B^{c})$ are measurable by assumption, the sets $B_{1},\ldots,B_{4}$ build a disjoint measurable partition of $E$. Based on this we may define a family $(T_{y})_{y \in E}$ of $\overbar{\mu}$-preserving transformations $T_{y}$ on $\overbar{X}$ by setting $T_{y} := \sigma$ if $y \in B_{2} \cup B_{4}$ and $T_{y} := \text{Id}$ if $y \in B_{1} \cup B_{3}$. By the ergodicity of $S$ either $B_{2}$ or $B_{4}$ must have positive measure (it is easily seen that otherwise $B \times \Omega$ would be a non-trivial $S$-invariant set). This implies that every measurable set $A \subseteq X$, which is invariant under the family $(T_{y})_{y \in E}$, is also invariant under $\sigma$ and thus trivial. In particular, the family is ergodic. 
	
	Let $T$ denote the step skew product arising from $S$ and $(T_{y})_{y \in E}$. We claim that $T$ is not ergodic. To show this we consider the sets $D_{1},\ldots,D_{4} \subseteq \Omega \times X$ given by
	\[D_{1} := (B_{1} \times \Omega) \times \{1\},~D_{2} := (B_{2} \times \Omega) \times \{1\},\]
	\[D_{3} := (B_{3} \times \Omega) \times \{2\},~D_{4} := (B_{4} \times \Omega) \times \{2\}\]
	and set $D := D_{1} \dot{\cup} \cdots \dot{\cup} D_{4}$. By the definition of the sets $B_{1},\ldots,B_{4}$ we obtain then
	\[D = \big((B \times \Omega) \times \{1\}\big) ~\dot{\cup}~ \big((B^{c} \times \Omega)\times \{2\}\big).\]
	We claim that $D$ is a $T$-invariant, non trivial set. To see this observe that 
	\begin{align*}
		\nu \otimes \overbar{\mu}(D) = m(B)\overbar{\mu}(1) + m(B^{c})\overbar{\mu}(2) = \frac{1}{2}.
	\end{align*}
	So it remains to show that $D$ is $T$-invariant. Let $D^{*}_{1}$ denote the set $(B_{1} \times B \times \Omega) \times \{1\}$. It is easy to see that $TD_{1}^{*} = (B \times \Omega) \times \{1\} \subseteq D$. Moreover, we have
	\begin{align*}
		\nu \otimes \overbar{\mu}\big(D_{1}^{*}\big) &= \nu(B_{1} \times B \times \Omega)~\overbar{\mu}(1) = \frac{1}{2}\int_{B_{1}}\pi(y,B)~dm(y) \\ &= \frac{1}{2}m(B_{1}) = \nu(B_{1} \times \Omega)~\overbar{\mu}(1) = \nu \otimes \overbar{\mu}(D_{1}).
	\end{align*}
	Since $D_{1}^{*} \subseteq D_{1}$ this shows that $\nu \otimes \overbar{\mu}(D_{1} \cap T^{-1}D^{c}) = 0$. By a similar argument we obtain $\nu \otimes \overbar{\mu}(D_{3} \cap T^{-1}D^{c}) = 0$. Let $D_{2}^{*}$ denote the set $(B_{2} \times B^{c} \times \Omega) \times \{1\}$. Then we have $TD_{2}^{*} = (B^{c} \times \Omega) \times \{2\} \subseteq D$. Moreover, we obtain
	\begin{align*}
		\nu \otimes \overbar{\mu}\big(D_{2}^{*}\big) &= \nu(B_{2} \times B^{c} \times \Omega)~\overbar{\mu}(1) = \frac{1}{2}\int_{B_{2}}\pi(y,B^{c})~dm(y) \\ &= \frac{1}{2}m(B_{2}) = \nu(B_{2} \times \Omega)~\overbar{\mu}(1) = \nu \otimes \overbar{\mu}(D_{2}).
	\end{align*}
	Since $D_{2}^{*} \subseteq D_{2}$ this shows that $\nu \otimes \overbar{\mu}(D_{2} \cap T^{-1}D^{c}) = 0$. By a similar argument we obtain $\nu \otimes \overbar{\mu}(D_{4} \cap T^{-1}D^{c}) = 0$. Noting that
	\[D \cap T^{-1}D^{c} = \bigcup_{i=1}^{4}(D_{i} \cap T^{-1}D^{c})\]
	we may conclude that $\nu \otimes \overbar{\mu}(D \cap T^{-1}D^{c}) = 0$. Since $T$ is measure preserving, this implies that $D$ is $T$-invariant. \hfill $\Box$ 
	
	\section{Random ergodic theorems} \label{SectionErgodicTheorem}
	
	In this final section we apply the results of the foregoing sections to obtain random ergodic theorems for measure preserving transformations, which are selected by Markov chains. The first result generalizes the classical result of Kakutani and Ryll-Nardzweski to this setting.
	
	\begin{theorem} \label{random}
		\thm Consider a stationary Markov chain $\xi = (\xi_{n})_{n=0}^{\infty}$ with state space $E$ arising from a Markov kernel $\pi$ and a $\pi$-invariant probability measure $m$. Consider furthermore a measurable family $(T_{y})_{y \in E}$ of measure preserving transformations on a probability space $(X,\mu)$. Then for every function $f \in L^{1}(\mu)$ there is a function $\smash{\widehat{f}} \in L^{1}(m \otimes \mu)$ such that almost surely the ergodic averages 
		\[A_{n}f := \frac{1}{n}\sum_{i=0}^{n-1}f \circ T_{\xi}^{~\!\!i}\]
		of $f$ along the random iterations $T_{\xi}^{i} := T_{\xi_{i-1}} \circ \cdots \circ T_{\xi_{0}}$ converge $\mu$-a.s. and in $L^{1}(\mu)$ to $\smash{\widehat{f}(\xi_{0},\cdot)}$. Moreover, if $\pi$ is strictly irreducible with respect to $m$, then for $m$-almost all $y \in E$ the function $\smash{\widehat{f}(y,\cdot)}$ coincides with the conditional expectation of $f$ under $\Sigma$.
	\end{theorem}
	
	\textit{Proof:} Let $(\Omega,\nu,S)$ denote the Markov shift corresponding to $\pi$ and $m$ and let $T$ be the step skew product arising from $S$ and $(T_{y})_{y \in E}$. For a given function $f \in L^{1}(\mu)$ we shall consider the functions $f^{*} := \1 \otimes f$ and $\smash{\overbar{f}} := \mathbb{E}[f^{*}|\Sigma_{T}]$, which are obviously in $L^{1}(\nu \otimes \mu)$. Then $\smash{\overbar{f}}$ is $T$-invariant, so by Proposition \ref{representation1} there is a function $\smash{\widehat{f}} \in L^{1}(m \otimes \mu)$ such that
	\[\overbar{f}(\omega,x) = \widehat{f}(\omega_{0},x)\]
	for $\nu \otimes \mu$-almost all ${(\omega,x) \in \Omega \times X}$. Writing
	\begin{align} \label{eq5}
		A_{n}f(\omega,x) = \frac{1}{n}\sum_{i=0}^{n-1}f \circ T_{\omega}^{~\!\!i}(x) = \frac{1}{n}\sum_{i=0}^{n-1}f^{*} \circ T^{i}(\omega,x)
	\end{align}
	we obtain thus
	\begin{align} \label{eq1}
		\lim_{n \to \infty}A_{n}f(\omega,x) =  \overbar{f}(\omega,x) = \widehat{f}(\omega_{0},x)
	\end{align}
	for $\nu \otimes \mu$-almost all $(\omega,x) \in \Omega \times X$ by Birkhoff's ergodic theorem. This implies the first convergence statement.
	
	For the mean convergence we shall first assume that $f \in L^{\infty}(\mu)$. Then obviously $f^{*} \in L^{\infty}(\mu)$ and $\|f^{*}\|_{\infty} \varleq \max\{1,\|f\|_{\infty}\}$. It is also not difficult to see that $\smash{\overbar{f}} \in L^{\infty}(\nu \otimes \mu)$ with $\|\smash{\overbar{f}}\|_{\infty} \varleq \|f^{*}\|_{\infty}$. Using these observations and the $T$-invariance of $\smash{\overbar{f}}$ we obtain
	\begin{align*}
		\big|A_{n}f(\omega,x) - \widehat{f}(\omega_{0},x)\big| &= \bigg|\frac{1}{n}\sum_{i=0}^{n-1}f^{*} \circ T^{i}(\omega,x) - \frac{1}{n}\sum_{i=0}^{n-1}\smash{\overbar{f}} \circ T^{i}(\omega,x)\bigg| \\ &\varleq \frac{1}{n}\sum_{i=0}^{n-1}\big|(f^{*} - \smash{\overbar{f}}) \circ T^{i}(\omega,x)\big| \varleq \|f^{*} - \smash{\overbar{f}}\|_{\infty} \\
		&\varleq 2\|f^{*}\|_{\infty} \varleq 2\max\{1,\|f\|_{\infty}\}
	\end{align*}
	for $\nu \otimes \mu$-almost all $(\omega,x) \in \Omega \times X$. Since the left hand side converges $\mu$-almost surely to zero for $\nu$-almost every $\omega \in \Omega$, we may apply dominated convergence to obtain
	\[\lim_{n \to \infty}\big\|A_{n}f(\omega,\cdot) - \widehat{f}(\omega_{0},\cdot)\big\|_{1} = 0\]
	for every such $\omega \in \Omega$. Now for a general $f \in L^{1}(\mu)$ a standard $3\varepsilon$-argument shows that for $\nu$-almost all $\omega \in \Omega$ the averages $A_{n}f(\omega,\cdot)$ form a Cauchy sequence and converge therefore to some limit in $L^{1}(\mu)$, which equals $\smash{\widehat{f}(\omega_{0},\cdot)}$ by (\ref{eq1}).
	
	Finally, if $\pi$ is strictly irreducible with respect to $m$, then $\pi$ is in particular irreducible with respect to $m$ by Proposition \ref{StrictlyIrreducibleImpliesIrreducible}. Accordingly, by Proposition \ref{ergodicity}, the shift $S$ is ergodic, thus $\Sigma_{S}$ is trivial. Using this we may apply statement $(ii)$ of Theorem \ref{main} to obtain that
	\begin{align*}
		\widehat{f}(\omega_{0},x) &= \mathbb{E}[f^{*}|\Sigma_{T}](\omega,x) = \mathbb{E}[\1 \otimes f|\Sigma_{S} \otimes \Sigma](\omega,x) \\ &= \mathbb{E}[\1|\Sigma_{S}](\omega)~\!\mathbb{E}[f|\Sigma](x) = \mathbb{E}[f|\Sigma](x)
	\end{align*}
	for $\nu \otimes \mu$-almost all $(\omega,x) \in \Omega \times X$, which shows the claim by the fact that $m$ is the push-forward measure of $\nu$ under $\xi_{0}$. \hfill $\Box$ \\

	We may also give an integrated version of the above random ergodic theorem. More precisely, instead of the time averaging operators $A_{n}$ we shall consider the expectation operators $M_{n}$ given by 
	\[M_{n}f(x) := \int_{\Omega}f \circ T_{\omega}^{~\!\!n}(x)~d\nu(\omega)\]
	for $n \in \mathbb{N}_{0}$, $f \in L^{1}(\mu)$ and $x \in X$. In case of a finite state space the operator $M_{n}$ describes just the weighted sum of $f$ over all possible compositions of transformations of length $n$. In this setting the theorem below was proven by Bufetov, cf. \cite[Theorem 3]{Buf00a} and \cite[Corollary 2]{Buf01}. Recall the definition of the space $L \log L(\mu)$, which consists of all measurable functions $f$ such that $\int |f|\log^{+}|f|~d\mu < \infty$.
	
	\begin{theorem} \label{random2}
		\thm Consider a stationary Markov chain $\xi = (\xi_{n})_{n=0}^{\infty}$ with state space $E$ arising from a Markov kernel $\pi$ and a $\pi$-invariant probability measure $m$. Consider furthermore a measurable family $(T_{y})_{y \in E}$ of measure preserving transformations on a probability space $(X,\mu)$. Then for every function $f \in L^{1}(\mu)$ there is a function $\smash{\widetilde{f}} \in L^{1}(\mu)$ such that we have
		\[\lim_{n \to \infty}\frac{1}{n}\sum_{j=0}^{n-1}M_{j}f = \widetilde{f}\]
		in $L^{1}(\mu)$. For $f \in L \log L(\mu)$ we obtain in addition $\mu$-a.s. convergence. Finally, if $\pi$ is strictly irreducible with respect to $m$, then the limit $\smash{\widetilde{f}}$ equals the conditional expectation of $f$ under $\Sigma$.
	\end{theorem}
	
	\textit{Proof:} As before let $(\Omega,\nu,S)$ denote the Markov shift corresponding to $\pi$ and $m$ and let $T$ be the step skew product arising from $S$ and $(T_{y})_{y \in E}$. Fix $f \in L^{1}(\mu)$. Then by Theorem \ref{random} we have $\nu \otimes \mu$-almost sure convergence of the averages ${A}_{n}f$ to a limit $\smash{\widehat{f}}$ in $L^{1}(m \otimes \mu)$. Setting $f^{*} := \1 \otimes f$ and $\smash{\overbar{f}} := \mathbb{E}[f^{*}|\Sigma_{T}]$ we obtain as above $\smash{\widehat{f}}(\omega_{0},x) = \overbar{f}(\omega,x)$ for $\nu \otimes \mu$-almost all $(\omega,x) \in \Omega \times X$. Note that by Fubini's theorem the partial integral 
	\[\widetilde{f}(\cdot) := \int_{E}\widehat{f}(y,\cdot)~dm(y) = \int_{\Omega}\overbar{f}(\omega,\cdot)~d\nu(\omega)\] 
	defines a function in $L^{1}(\mu)$. (The last equality follows again from the fact that $m$ is the push-forward of $\nu$ under the projection $\xi_{0}$.) Thus we obtain
	\begingroup
	\allowdisplaybreaks
	\begin{align*}
		\bigg\|\frac{1}{n}\sum_{j=0}^{n-1}M_{j}f - \widetilde{f}~\bigg\|_{1} &= \int_{X}\bigg|\int_{\Omega}\bigg(\frac{1}{n}\sum_{j=0}^{n-1}f \circ T_{\omega}^{j}(x) - \overbar{f}(\omega,x)\bigg)~d\nu(\omega)\bigg|~d\mu(x) \\
		&\varleq \int_{\Omega}\int_{X}\bigg|\frac{1}{n}\sum_{j=0}^{n-1}f \circ T_{\omega}^{j}(x) - \overbar{f}(\omega,x)\bigg|~d\mu(x)d\nu(\omega) \\
		&= \int_{\Omega \times X}\bigg|\frac{1}{n}\sum_{j=0}^{n-1}f^{*} \circ T^{j}(\omega,x) - \overbar{f}(\omega,x)\bigg|~d\nu \otimes \mu(\omega,x) \\
		&= \bigg\|\frac{1}{n}\sum_{j=0}^{n-1}f^{*} \circ T^{j} - \overbar{f}\bigg\|_{1}~,
	\end{align*}
	\endgroup
	where the latter term converges to zero by the mean ergodic theorem. If $\pi$ is strictly irreducible with respect to $m$, we have furthermore
	\[\widetilde{f}(x) = \int_{E}\widehat{f}(y,x)~dm(y) = \int_{E}\mathbb{E}[f|\Sigma](x)~dm(y) = \mathbb{E}[f|\Sigma](x)\]
	for $\mu$-almost all $x \in X$ by Theorem \ref{random}. 
	
	Now assume $f \in L \log L(\mu)$. Then $f$ is integrable with respect to $\mu$, which in turn gives that $f^{*}$ and therefore $A_{n}f$ is integrable with respect to $\nu \otimes \mu$ for all $n \in \mathbb{N}$. Again by Fubini's theorem we may conclude that $A_{n}f(\cdot,x)$ is a function in $L^{1}(\nu)$ satisfying
	\begin{align} \label{eq4}
		\int_{\Omega} A_{n}f(\omega,x)~d\nu(\omega) = \frac{1}{n}\sum_{j=0}^{n-1}\int_{\Omega} f \circ T_{\omega}^{j}(x)~d\nu(\omega) =\frac{1}{n}\sum_{j=0}^{n-1}M_{j}f(x)
	\end{align}
	for $\mu$-almost all $x \in X$. Consider the maximal function $F\colon \Omega \times X \to [0,\infty)$ given by
	\[F(\omega,x) := \sup_{n \in \mathbb{N}} |A_{n}f(\omega,x)|\]
	for $(\omega,x) \in \Omega \times X$. Since $f \in L \log L(\mu)$, we have $f^{*} \in L \log L(\nu \otimes \mu)$. Therefore, by (\ref{eq5}) and Wiener's dominated ergodic theorem, cf. \cite[Chapter 1, Theorem 6.3]{Kre85}, we obtain $F \in  L^{1}(\nu \otimes \mu)$. This implies that $F(\cdot,x)$ is in $L^{1}(\nu)$ for $\mu$-almost all $x \in X$. In particular, for $\mu$-almost all $x \in X$ the functions $A_{n}f(\cdot,x)$ converge $\nu$-almost surely to $\overbar{f}(\cdot,x)$ and is dominated by the integrable function $F(\cdot,x)$. Dominated convergence thus yields
	\[\lim_{n \to \infty}\frac{1}{n}\sum_{j=0}^{n-1}M_{j}f(x) =
	\int_{\Omega}\overbar{f}(\omega,x)~d\nu(\omega) = \widetilde{f}(x)\]
	for $\mu$-almost all $x \in X$ by (\ref{eq4}). \hfill $\Box$ \\ 
	
	\medskip
	
	\textbf{Outlook.} It is well known that in the case of an i.i.d. selection of transformations the assumptions of the random ergodic theorem may be weakened. In fact, instead of a common invariant probability measure it suffices to consider a probability measure, which is stationary for the family of transformations, cf. \cite[Chapter 1.2]{Kif86}. The latter means that the measure is invariant under the Markov kernel induced by the random transformations on the state space of the random dynamical system. In fact, one can show that a measure is stationary if and only if the respective product measure is invariant under the skew product.  This observation is particular useful in situations where the transformations are not invertible and the assumption of a common invariant measure excludes interesting bevahior (as it is often the case in low dimensional dynamics). However, if we allow the selection process to be Markov chain, a similarly convenient correspondence is not at hand, so it seems reasonable to confine to invariant probability measures. As discussed in the introduction in many applications this assumption is naturally instantiated. However, in certain contexts it is more or less restrictive, see e.\@ g. \cite{Mat22}. It would therefore be interesting, to what extend this assumption can be relaxed also in the present setting. 
	
	\medskip
	
	In the case of a finite state Markov chain Theorem \ref{random2} can be improved in the sense that pointwise convergence holds for all $L^{1}$-functions, cf. \cite[Corollary 2]{Buf01}. This is due to the fact that the averages over all paths of the Markov chain with a fixed first symbol can be encoded by a convenient Markov operator (which is also implicitly present in the computations of the proof of Theorem \ref{representation2}). Since in this case the averages considered in Theorem \ref{random2} are just a finite weighted sum of the ergodic averages of this operator, one may use the general ergodic theorem for Dunford-Schwartz operators to obtain pointwise convergence for all $L^{1}$-functions. In the general setting this approach seems to fail. A priori there is no argument at hand allowing to interchange integrating and taking the limit of the operator ergodic averages. It would therefore be interesting if there is an alternative argument showing pointwise convergence for $L^{1}$-functions. 
	
	\medskip

	\textbf{Acknowledgements.} The authors would like to thank Alexander I. Bufetov and Anush Tserunyan for stimulating discussions, valuable hints and helpful remarks as well as the anonymous referee for useful suggestions and comments. Felix Pogorzelski und Elias Zimmermann gratefully acknowledge support by the German Israeli Foundation for Scientific Research and Development (GIF, grant I-1485-304.6/2019).

	\textsc{Mathematical Institute, University of Leipzig} 
	
	\textsc{Augustusplatz 10, 04109 Leipzig}
	
	\textit{p.lummerzheim@math.leidenuniv.nl}
	
	\textit{felix.pogorzelski@math.uni-leipzig.de}
	
	\textit{elias.zimmermann@math.uni-leipzig.de}

\end{document}